\documentclass[12pt]{article}

\newcommand{\Prob}{\mbox{\sc P}}
\newtheorem{thm}{Theorem}
\newtheorem{cor}{Corollary}

\begin{document}

\begin{center}
{\bf\LARGE  On Reichenbach's causal betweenness}\\
\vspace{1cm}
Va\v{s}ek Chv\'{a}tal\\
{\em
Canada Research Chair in Combinatorial Optimization\\
Department of Computer Science and Software Engineering\\
Concordia University\\
Montreal, Quebec H3G 1M8, Canada\\
}
\vspace{0.5cm}
Baoyindureng Wu\\
{\em
College of Mathematics and System Sciences\\
Xinjiang University\\
Urumqi 830046, PR China
}
\end{center}

\vspace{0.7cm}

\begin{center}
  {\bf Abstract}\\
  \vspace{0.3cm} {\small We characterize, by easily verifiable
    properties, abstract ternary relations isomorphic to the causal
    betweenness introduced by Hans Reichenbach.}
\end{center}

\vspace{0.3cm}

\section{Introduction}

A {\em finite probability space\/} is an ordered pair $(S, p)$ where
$S$ is a finite set and $p:S\rightarrow [0,1]$ is a function such that
\[
{\textstyle \sum_{s\in S}p(s)=1}.
\]
The set $S$ is called a {\em sample space\/} and its subsets are
called {\em events\/}; when $A$ and $B$ are events, $AB$ denotes
$A\cap B$. The {\em probability $\Prob(A)$ of event $A$\/} is defined
by
\[
{\textstyle \Prob(A)=\sum_{s\in A}p(s)};
\]
when $A$ and $B$ are events with $\Prob(B)>0$, the {\em conditional
  probability $\Prob(A|B)$ of $A$ given $B$\/} is defined by
\[
\Prob(A|B) = \frac{\Prob(AB)}{\Prob(B)}.
\]

\bigskip 

Hans Reichenbach (1956, p.~190) defined an event $B$ to be {\em
  causally between\/} events $A$ and $C$ if the following relations
hold:
\begin{eqnarray}
\Prob(AC)&>&\Prob(A)\cdot\Prob(C),
\label{cb1}\\
\rule{0pt}{12pt} \Prob(C|B)&>& \Prob(C|A),\hspace{0.4in}
\label{cb2}\\
\rule{0pt}{12pt} \Prob(A|B)&>& \Prob(A|C)
,\label{cb3}\\
\Prob(AC|B)&=&\Prob(A|B)\cdot\Prob(C|B)
,\label{cb4}
\end{eqnarray}
and 
\begin{equation}\label{cb5}
\Prob(B-A)>0,\;\Prob(B-C)>0.
\end{equation}
(The conditional probabilities in (\ref{cb2}) -- (\ref{cb4}) are
well-defined: (\ref{cb1}) guarantees that $\Prob(A)>0$, $\Prob(C)>0$
and (\ref{cb5}) guarantees that $\Prob(B)>0$.  In terms of
Reichenbach (1956, p.~189, pp. 201 -- 205) equation (\ref{cb4}) means
that $B$ {\em screens off\/} $A$ from $C$.)  Following Reichenbach's
work, causal betweenness was considered by von Bretzel (1977), Ellett
and
Ericson (1986), Dowe (1992), Weber (1997), Korb (1999), and others.\\

Given a set $X$ of events in a finite probability space, we let ${\cal
  CB}(X)$ denote the set of all ordered triples $(A,B,C)$ such that
$A,B,C$ are events in $X$ with properties (\ref{cb1}) -- (\ref{cb5}).
We say that a ternary relation $\cal B$ on a finite ground set is an
{\em abstract causal betweenness\/} if, and only if, there is a set
$X$ of events in a finite probability space such that ${\cal CB}(X)$
is isomorphic to $\cal B$. Our Theorem~\ref{thm1} characterizes
abstract causal betweennesses by easily verifiable properties.\\

We call a ternary relation $\cal B$ a {\em betweenness\/} if
\begin{eqnarray*}
\rule{0pt}{12pt} (A,B,C)\in{\cal B} &\Rightarrow& \mbox{$A,B,C$ are all distinct}
,\label{cb31}\\
\rule{0pt}{12pt} (A,B,C)\in{\cal B} &\Rightarrow& (C,B,A)\in{\cal B} 
,\label{cb32}\\
\rule{0pt}{12pt} (A,B,C)\in{\cal B} &\Rightarrow& (C,A,B)\not\in{\cal B}
.\label{cb33}
\end{eqnarray*}
The familiar concept of betweenness in Euclidean geometry generalizes
in diverse branches of mathematics to betweennesses ${\cal B}$ with
the property
\begin{equation}\label{trans}
(ABC), (ADB) \in {\cal B} \;\Rightarrow\; (ADC)\in {\cal B}.  
\end{equation} 
These relations include {\em metric betweennes\/}, {\em lattice
  betweenness\/} in modular lattices, and {\em algebraic betweenness}
(see, for instance, Pitcher and Smiley (1942), Smiley (1943),
Hashimoto (1958), Bumcrot (1964)).  Our Corollary~\ref{cor1} asserts
that every betweenness ${\cal B}$ with property (\ref{trans}) is an
abstract causal betweenness.\\

In Reichenbach's investigations, events occur in time and time order
is reduced to causal order (Reichenbach 1956, p.~24). An event that is
causally between events $A$ and $C$ does not necessarily occur between
$A$ and $C$ (after $A$ and before $C$, or else after $C$ and before
$A$). To elaborate on this point, we say that a ternary relation $\cal
B$ on a finite set is {\em totally orderable\/} if, and only if, there
is a mapping $t$ from the ground set of $\cal B$ to a set with a total
order $\prec$ such that
\[
(A,B,C)\in {\cal B} \;\;\Rightarrow\;\;
(t(A)\prec t(B)\prec t(C) \;\mbox{ or }\; t(C)\prec t(B)\prec t(A)).  
\]
Reichenbach (1956, p.~192) pointed out that $\cal B$ is not totally
orderable when it includes the ordered triples $(A_1,A_2,A_3)$,
$(A_1,A_2,A_4)$, $(A_4,A_2,A_3)$; these three triples, along with
their reversals $(A_3,A_2,A_1)$, $(A_4,A_2,A_1)$, $(A_3,A_2,A_4)$,
constitute an abstract causal betweenness. Therefore not every
abstract causal betweenness is totally orderable.\\

Opatrn\'{y} (1979) proved that recognizing totally orderable ternary
relations is hard: the problem is ${\cal NP}$-complete. (Readers
unfamiliar with the notion of ${\cal NP}$-completeness are referred to
the monograph of Garey \& Johnson (1979).) The problem does not get
any easier when its input is restricted to abstract causal
betweennesses: our Corollary~\ref{cor2} asserts that every totally
orderable betweenness is an abstract causal betweenness, and so
testing an arbitrary ternary relation $\cal B$ for total orderability
reduces to testing an abstract causal betweenness for total
orderability. (We first test $\cal B$ for being an abstract causal
betweenness; Theorem~\ref{thm1} shows how to carry out this test
easily; if $\cal B$ fails it, then Corollary~\ref{cor2} guarantees
that $\cal B$ is not totally orderable.)

\section{Results}\label{sec2}

With each betweenness ${\cal B}$ on a ground set $X$, we associate a
directed graph $G({\cal B})$. Its vertices are all two-point subsets
of $X$; its edges are all ordered pairs $(\{A,B\},\{A,C\})$ such that
$(ABC)\in {\cal B}$. These graphs may contain directed cycles: if
$\{(DAB), (DBC), (DCA)\}\subseteq {\cal B}$, then $G({\cal B})$ contains the
directed cycle
\[
\{D,A\} \rightarrow \{D,B\} \rightarrow \{D,C\} \rightarrow \{D,A\},
\]
if $\{(CAB), (DBC), (ACD), (BDA)\}\subseteq {\cal B}$, then $G({\cal B})$
contains the directed cycle
\[
\{A,B\} \rightarrow \{B,C\} \rightarrow \{C,D\} \rightarrow \{D,A\} \rightarrow \{A,B\},
\]
and so on. 

\begin{thm}\label{thm1} 
  A ternary relation $\cal B$ on a finite set is an abstract causal
  betweenness if and only if $\cal B$ is a betweenness and $G({\cal
    B})$ contains no directed cycle.
\end{thm}

\begin{cor}\label{cor1} 
  Every betweenness $\cal B$ with the property
\[
(ABC), (ADB) \in {\cal B} \;\Rightarrow\; (ADC)\in {\cal B}.  
\]
is an abstract causal betweenness.
\end{cor}

\begin{cor}\label{cor2} 
Every totally orderable betweenness is an abstract causal betweenness.
\end{cor}

\section{Proofs}

\noindent {\bf Proof of Theorem~\ref{thm1}.}

\rule{0pt}{20pt}
The ``if'' part: Consider an arbitrary betweenness ${\cal B}$ on a ground set $X$ such
that the directed graph $G({\cal B})$ contains no directed
cycle. Without loss of generality, $X=\{1,2,\ldots ,m\}$ for some
positive integer $m$.  We shall construct a finite probability space
and events $E_1,E_2,\ldots ,E_m$ in this space in such a way that
$E_j$ is causally between $E_i$ and $E_k$ if and only if
$(i,j,k)\in {\cal B}$.\\

The construction proceeds in two stages. First, we choose an
arbitrarily small positive $\varepsilon$ and we construct functions
\begin{eqnarray*}
\beta:  \{W: W\subseteq X, |W|=2\} &\rightarrow & (0.25,\;0.25+\varepsilon),\\
\gamma: \{W: W\subseteq X, |W|=3\} &\rightarrow & (0.125,\;0.125+\varepsilon)
\end{eqnarray*}
such that $(i,j,k)\in {\cal B}$ if and only if
\begin{eqnarray*}
\gamma(\{i,j,k\})&=&2\beta(\{i,j\})\beta(\{j,k\}),
\\
\beta(\{i,j\})&>& \beta(\{i,k\}),
\\
\beta(\{j,k\})&>& \beta(\{i,k\}).
\end{eqnarray*}
Then we construct a finite probability space and events
$E_1,E_2,\ldots ,E_m$ in this space in such a way that
\[
\begin{array}{rcll}
\Prob(E_i) &=& 0.5 &\mbox{ for all subscripts $i$,}\\
\Prob(E_iE_j) &=& \beta(\{i,j\}) &\mbox{ for all choices of distinct subscripts $i,j$,}\\
\Prob(E_iE_jE_k) &=& \gamma(\{i,j,k\}) &\mbox{ for all choices of distinct subscripts $i,j,k$.}\\
\end{array}
\]

\bigskip

In the first stage, we choose $\varepsilon$ and $\delta$ so that
\[
0<\varepsilon<m^{-2}4^{-m}, \;\;\; 0<\delta<m^{-2}\varepsilon.
\]
Since $G({\cal B})$ contains no directed cycle, there is a mapping
\[
\rho:\{W\subset X: |W|=2\}\rightarrow \{1,2,\ldots ,m(m-1)/2\}
\]
such that
\[
(i,j,k)\in {\cal B} \;\;\Rightarrow\;\; \rho(\{i,j\}) > \rho (\{i,k\}).
\]
We set
\[
\beta(\{i,j\})=0.25+\delta\rho(\{i,j\})
\]
for all choices of distinct subscripts $i,j$.  For all $(i,j,k)$ in
${\cal B}$, we set
\[
\gamma(\{i,j,k\})=2\beta(\{i,j\})\beta(\{j,k\});
\]
if $i,j,k$ are distinct subscripts such that none of $(i,j,k)$,
$(j,k,i)$, $(k,i,j)$ is in ${\cal B}$, then we choose
$\gamma(\{i,j,k\})$ in the interval $(0.125,\;0.125+\varepsilon)$ and
distinct from all three of
\[
2\beta(\{i,j\})\beta(\{j,k\}), \;\;
2\beta(\{j,k\})\beta(\{k,i\}), \;\;
2\beta(\{k,i\})\beta(\{i,j\}). 
\]
The upper bound on $\delta$ guarantees that
\[
0.25 < \beta(\{i,j\}) <0.25 +\varepsilon/2
\]
for all choices of distinct subscripts $i,j$ and that
\[
0.125 < \gamma(\{i,j,k\}) < 0.125 +\varepsilon
\]
for all choices of distinct subscripts $i,j,k$.\\

In the second stage, we begin with sample space $\{0,1\}^m$ and we set 
\[
E_i=\{s\in \{0,1\}^m: s_i=1\}\;\;\; (i=1,2,\ldots ,m).
\]
For each subset $W$ of $X$, let $\chi^W$ denote the vector in
$\{0,1\}^m$, whose $i$-th coordinate is $1$ if and only if $i\in W$.
We will complete the proof by exhibiting a function $p:\{0,1\}^m\rightarrow
(0,1)$ such that 
\begin{equation}\label{3}
{\textstyle \sum(p(\chi^W):i,j,k\in W) = \gamma(\{i,j,k\})}
\end{equation}
for all choices of distinct subscripts $i,j,k$, 
\begin{equation}\label{2}
{\textstyle \sum(p(\chi^W):i,j\in W) = \beta(\{i,j\})}
\end{equation}
for all choices of distinct subscripts $i,j$,
\begin{equation}\label{1}
{\textstyle \sum(p(\chi^W):i\in W) = 0.5}
\end{equation}
for all subscripts $i$, and 
\begin{equation}\label{0}
{\textstyle \sum_W p(\chi^W)=1.}
\end{equation}
For this purpose, we set first
\[
p(\chi^W)=2^{-m}\;\;\;\mbox{ whenever $|W|\ge 4$}
\]
and
\[
p(\chi^{\{i,j,k\}})=2^{-m}+(\gamma(\{i,j,k\})-0.125)
\]
to satisfy (\ref{3}) for all choices of distinct subscripts $i,j,k$, then
\begin{eqnarray*}
p(\chi^{\{i,j\}}) &=& 2^{-m} + (\beta(\{i,j\})-0.25)\\
&& - \;{\textstyle\sum(\,p(\chi^W)-2^{-m}\!:\; W\supset \{i,j\},\, |W|=3)}
\end{eqnarray*}
to satisfy (\ref{2}) for all choices of distinct subscripts $i,j$, then 
\[
p(\chi^{\{i\}}) = 2^{-m} -{\textstyle\sum(\,p(\chi^W)-2^{-m}\!:\; W\ni i,\, 2\le |W|\le 3)}
\]
to satisfy (\ref{1}) for all subscripts $i$, and finally
\[
p(\chi^\emptyset)=1-{\textstyle \sum_{W\ne \emptyset} p(\chi^W).}
\]
to satisfy (\ref{0}). Now
\[
2^{-m} \;<\; p(\chi^{\{i,j,k\}}) \;<\; 2^{-m}+\varepsilon
\]
for all choices of distinct subscripts $i,j,k$, 
\[
2^{-m} -m\varepsilon \;<\; p(\chi^{\{i,j\}}) \;<\; 2^{-m}+\varepsilon
\]
for all choices of distinct subscripts $i,j$, 
\[
2^{-m} -m^2\varepsilon \;<\; p(\chi^{\{i\}}) \;<\; 2^{-m}+ m^2\varepsilon
\]
for all subscripts $i$, and so the upper bound on $\varepsilon$ guarantees that
all $p(\chi^W)$ are positive.\\

The ``only if'' part: Consider an arbitrary set $X$ of events in a
finite probability space. Reichenbach (1956, p.~191) proved that
${\cal CB}(X)$ is a betweenness. We shall reproduce his argument here
and we shall show that $G({\cal CB}(X))$ contains no directed cycle.\\

To prove that ${\cal CB}(X)$ is a betweenness, consider a triple
$(A,B,C)$ of events that satify (\ref{cb1}) -- (\ref{cb5}). Assumption
(\ref{cb5}) guarantees $B\ne A$ and $B\ne C$; assumption (\ref{cb2})
guarantees $\Prob(C|A)<1$, which implies $C\ne A$; now $A,B,C$ are all
distinct. Since the set of assumptions (\ref{cb1}) -- (\ref{cb5}) is
invariant under the switch $A\leftrightarrow C$, the triple $(C,B,A)$
satisfies them in place of $(A,B,C)$. Finally, the triple $(C,A,B)$
fails to satisfy (\ref{cb4}) in place of $(A,B,C)$ since (\ref{cb4})
and (\ref{cb2}) imply
\[
\Prob(CB|A) \;=\; \Prob(C|B)\cdot\Prob(B|A) \;>\; \Prob(C|A)\cdot\Prob(B|A).
\]

\bigskip

To see that $G({\cal CB}(X)$ contains no directed cycle, observe that
(\ref{cb3}) implies
\[
(A,B,C)\in {\cal B} \;\;\Rightarrow\;\; 
\frac{\Prob(AB)}{\Prob(A)\cdot\Prob(B)} \;>\; \frac{\Prob(AC)}{\Prob(A)\cdot\Prob(C)}.
\]

\noindent {\bf Proof of Corollary~\ref{cor1}.}
Writing $\sigma(A,B)$ for the number of elements $D$ such that $(ADB)
\in {\cal B}$, observe that (\ref{trans}) implies 
\[
(A,B,C)\in {\cal B} \;\;\Rightarrow\;\; \sigma(A,B) < \sigma (A,C),
\]
and so $G({\cal B})$ contains no directed cycle.\\ 

\noindent {\bf Proof of Corollary~\ref{cor2}.}
Consider an arbitrary betweenness $\cal B$ on a finite set along with
a mapping $t$ from the ground set of $\cal B$ to the set of real
numbers such that
\[
(A,B,C)\in {\cal B} \;\;\Rightarrow\;\;
(t(A)< t(B)< t(C) \;\mbox{ or }\; t(C)< t(B)< t(A)).  
\]
Writing $\tau(A,B)=|t(A)-t(B)|$, observe that
\[
(A,B,C)\in {\cal B} \;\;\Rightarrow\;\; \tau(A,B) < \tau (A,C),
\]
and so $G({\cal B})$ contains no directed cycle. 

\bigskip

\begin{center}
{\bf Acknowledgment}
\end{center}

\noindent 
This note was written while the second author was visiting the first
author in ConCoCO (Concordia Computational Combinatorial Optimization
Laboratory). The research was undertaken, in part, thanks to funding
from the Canada Research Chairs Program and by the China Scholarship Council.\\

\bigskip

\noindent 
{\bf\Large  References}\\

\noindent 
von Bretzel, P. (1977). Concerning a probabilistic theory of
causation adequate for the causal theory of time. {\em Synthese, 35,\/} 173 -- 190.\\

\noindent 
Bumcrot, R.J. (1964). Betweenness geometry in lattices.  {\em
  Rendiconti del Circolo Matematico di Palermo, 13,\/} 11 -- 28.\\

\noindent 
Dowe P. (1992). Process causality and asymmetry. {\em
  Erkenntnis, 37,\/} 179 -- 196.\\

\noindent 
Ellett, Jr., F.S. \& Ericson, D.P. (1986). Correlation, partial
correlation, and causation. {\em Synthese, 67,\/} 157 -- 173.\\

\noindent 
Garey, M.R. \& Johnson, D.S. (1979). {\em Computers and
  Intractability: A Guide to the Theory of NP-Completeness.\/}
(San Francisco: W.H. Freeman and Company.)\\

\noindent 
Hashimoto, J. (1958). Betweenness geometry. {\em Osaka Mathematical Journal, 10,\/} 147 -- 158.\\

\noindent 
Korb, K.B. (1999). Probabilistic causal structure. (In H. Sankey
(Ed.), {\em Causation and Laws of Nature\/} (pp. 265 -- 311).  The
Netherlands: Kluwer Academic Publishers.)\\

\noindent 
Opatrn\'{y}, J. (1979). Total ordering problem.
{\em SIAM Journal on Computing, 8,\/} 111 -- 114.\\

\noindent 
Pitcher, E. \& Smiley, M.F. (1942). Transitivities of
betweenness. {\em Transactions of the American Mathematical Society, 52,\/} 95  -- 114.\\

\noindent 
Reichenbach, H. (1956). {\em The Direction of Time.\/} (Berkeley and Los Angeles: University of California Press.)\\

\noindent 
Smiley, M.F. (1943). A comparison of algebraic, metric, and lattice
betweenness. {\em Bulletin of the American Mathematical Society, 49,\/} 246 -- 252.\\

\noindent 
Weber, G.D. (1997). Discovering causal relations by experimentation: Causal trees. 
(In {\em Proceedings of the Eighth Midwest AI and Cognitive Science Conference\/} (pp. 91 -- 98). 
Menlo Park: Association for the Advancement of Artificial Intelligence.)\\

\end{document}